\documentclass[11pt]{article}
\usepackage{amsfonts}

\usepackage{graphics}
\usepackage{indentfirst}
\usepackage{cite}
\usepackage{latexsym}
\usepackage{amsmath}
\usepackage{amssymb}
\usepackage[dvips]{epsfig}
\usepackage{amscd}

\newtheorem{theorem}{Theorem}[section]
\newtheorem{remark}{Remark}[section]

\newtheorem{definition}{Definition}[section]
\newtheorem{lemma}[theorem]{Lemma}

\newtheorem{proposition}[theorem]{Proposition}

\newcommand{\n}{\rho}
\newcommand{\rt}{R_T}
\newcommand{\ti}{\tilde}

\newcommand{\lm}{\lambda}
\newcommand\divg{{\text{div}}}
\def\pf{{\it Proof.}  }
\renewcommand{\div}{ {\rm div }  }

\newcommand{\na}{\nabla }
\newcommand{\vp}{\varphi }

\newcommand{\pa}{\partial}

\newcommand{\bi}{\bibitem}

\newcommand{\bt}{\begin{theorem}}
\newcommand{\bl}{\begin{lemma}}
\newcommand{\el}{\end{lemma}}
\newcommand{\et}{\end{theorem}}
\newcommand{\ga}{\gamma}

\newcommand{\al}{a }
\newcommand{\de}{\delta}
\newcommand{\ve}{\varepsilon}
\newcommand{\la}{\label}
\newcommand{\si}{\sigma}

\newcommand{\bn}{\begin{eqnarray}}
\newcommand{\en}{\end{eqnarray}}
\newcommand{\bnn}{\begin{eqnarray*}}
\newcommand{\enn}{\end{eqnarray*}}

\newcommand{\bnnn}{\begin{eqnarray*}}
\newcommand{\ennn}{\end{eqnarray*}}
\newcommand{\ben}{\begin{enumerate}}
\newcommand{\een}{\end{enumerate}}

\newcommand{\ba}{\begin{aligned}}
\newcommand{\ea}{\end{aligned}}
\newcommand{\be}{\begin{equation}}
\newcommand{\ee}{\end{equation}}

\def\O{\mathbb{R}^2}
\def\p{\partial}
\def\norm[#1]#2{\|#2\|_{#1}}

\def\lap{\triangle}
\def\g{\gamma}
\def\lam{\lambda}

\def\o{\omega}
\def\rr{\mathbb{R}^2}

\makeatletter      
\@addtoreset{equation}{section}
\makeatother       

\title{ Global Classical Solution to the Cauchy Problem of 2D Baratropic Compressible Navier-Stokes System with Vacuum and Large Initial Data  \thanks{The research of X. Huang was partially supported by NNSFC 11101392.
 Email: xdhuang@amss.ac.cn (X. Huang), ajingli@gmail.com (J. Li).
 }}

\date{}
\author{Xiangdi H{\small UANG}$^{a}$,  Jing L{\small I}$^{b}$  \\[3mm] {\normalsize $^a$ Academy of Mathematics and System Sciences,} \\
{\normalsize Chinese Academy of Sciences, Beijing 100190, P. R. China} \\[2mm]
{\normalsize $^b$ Institute of Applied Mathematics, AMSS,} \\ {\normalsize \&   Hua Loo-Keng Key Laboratory of Mathematics,}\\
{\normalsize  Chinese Academy of Sciences,    Beijing 100190,
P. R. China}
 }

\begin{document}
\maketitle

\begin{abstract}
 The authors establish the global existence and uniqueness of strong and classical solutions to the Cauchy problem for the barotropic compressible
Navier-Stokes equations on the whole  two-dimensional space with  vacuum  as far field density and
    with    no restrictions on the size of   initial data  provided the shear  viscosity   is a positive constant and the  bulk one is $\lambda = \rho^{\beta}$ with $\beta>4/3$.
\end{abstract}


\section{Introduction and main results}
We are concerned with   the two-dimensional barotropic  compressible Navier-Stokes equations which read as follows:
\be\la{n1}
\begin{cases} \rho_t + \div(\rho u) = 0,\\
 (\rho u)_t + \div(\rho u\otimes u) + \nabla P = \mu\lap u + \nabla((\mu + \lam)\div u),
\end{cases}
\ee
where   $t\ge 0, x=(x_1,x_2)\in \Omega\subset\O, \rho=\n(x,t)$ and $u=(u_1(x,t),u_2(x,t))$ represent, respectively, the density and velocity, and the pressure $P$ is given by
\be\la{n2}
P(\rho) = R\rho^{\gamma},\quad \ga>1.
\ee
The shear viscosity   $\mu$ and the bulk one $\lambda$ satisfy  the following hypothesis:
\be\la{n3}
0<\mu = const,\quad \lam(\n)  = b\rho^{\beta},\,\, b>0,\,\, \beta>0.
\ee
 In the sequel, we set $R=b = 1$ without loss of any generality.
 Let $\Omega=\O$ and we consider the Cauchy problem with $(\n,u)$ vanishing at infinity (in some weak sense). For given initial data $\n_0$ and $u_0,$   we require that
\be \la{n4} \n(x,0)=\n_0(x), \quad \n u(x,0)=\n_0u_0(x),\quad x\in \Omega= \O.\ee

When both the shear and bulk viscosities are positive constants, there is a huge literature  on the global existence and large-time
behavior of solutions to (\ref{n1}). The one-dimensional problem
has been studied extensively by many people, see
\cite{Kaz,Ser1,Ser2,Hof} and the references therein.  For the
multi-dimensional case, when the initial data $\n_0,m_0$ are sufficiently regular and the initial density $\n_0$ has a positive lower bound, the local existence and uniqueness of
classical solutions are known in \cite{Na,se1}.  Recently, for the case that the initial density need not be positive
and may vanish in open sets,   the existence and uniqueness of local strong and classical solutions were obtained by \cite{cho1,
K2, sal}.  More recently,   for two-dimensional case and $\Omega=\O,$ Li-Liang \cite{hlma} obtained the existence and uniqueness of the  local strong and classical solutions to \eqref{n1}-\eqref{n4} with vacuum  as far field density.
The global classical solutions were
first obtained by Matsumura-Nishida \cite{M1} for initial data
close to a non-vacuum equilibrium in some Sobolev space $H^s.$   Such theory was later generalized to weak solutions by  Hoff \cite{Ho4}
 and solutions in Besov
spaces by Danchin \cite{da1}.   For the
existence of solutions for large data,  the major breakthrough is due to
Lions \cite{L1} (see also  Feireisl \cite{F1,Fe}), where he obtained
global existence of weak solutions  when the exponent $\ga$ is suitably large. The
main restriction on initial data is that the initial energy is
finite, so that the   density is allowed to vanish initially.  Recently, Huang-Li-Xin  \cite{hlx1} established the global existence and uniqueness of classical
solutions to the Cauchy problem for the isentropic compressible Navier-Stokes equations in three-dimensional space with smooth initial data which are of small energy  but possibly large oscillations; in particular, the initial density is allowed to vanish, even has compact support.

  However,  there are few results regarding global  strong  solvability for equations of multi-dimensional motions of viscous gas  with    no restrictions on the size of   initial data. One of the first ever ones is due to  Vaigant-Kazhikhov \cite{Ka} who obtained   that   the two-dimensional system  \eqref{n1}-\eqref{n4} admits a unique global strong solution for large  initial data away from vacuum provided $\beta>3$ and the domain $\Omega$ is bounded. Recently,   under some additional compatibility  conditions on the periodic initial data,  Jiu-Wang-Xin \cite{jwx} considered periodic classical solutions and  removed  the  condition that the initial density should be away from vacuum in  Vaigant-Kazhikhov \cite{Ka} but still  under the  same condition  $\beta>3$ as that in \cite{Ka}. More recently,
  for  periodic   initial data with initial density allowed to vanish,  we   \cite{hlia} not only relax the crucial condition $\beta>3$ of  \cite{Ka} to the one that  $\beta>4/3,$ but also obtain both the time-independent  upper bound of the density and the large-time   behavior of the strong and weak solutions.
It should be noted that    \cite{Ka,jwx,hlia} only consider the periodic case  or the case of  bounded domains and the global existence of strong and classical solutions to the Cauchy problem \eqref{n1}-\eqref{n4} in the whole space $\O$ remains open. In fact, this is  the aim of this paper.

Before stating the main results, we explain the notations and
conventions used throughout this paper. We denote
\bnn\int fdx=\int_{\O}fdx  .\enn For $1\le r\le \infty ,$  we also denote the standard Lebesgue and
  Sobolev spaces as follows:
$$ L^r=L^r(\O), \quad W^{s,r}= W^{s,r}(\O),  \quad H^s= W^{s,2} .  $$

Then, we give the definition of strong solutions to \eqref{n1}:
\begin{definition} If  all derivatives involved in \eqref{n1} for $(\rho,u)  $  are regular distributions, and   equations  \eqref{n1} hold   almost everywhere   in $\rr\times (0,T),$ then $(\n,u)$  is called a  strong solution to  \eqref{n1}.
\end{definition}
Thus, the first main result concerning the global existence   of  strong solutions can be stated as follows:
\begin{theorem}\la{t2} Assume that \be\la{bet}\beta>4/3,\quad\ga>1,\ee
 and that the initial data $(0\le \n_0,u_0)$ satisfy that for some $q>2$ and $a\in (1,2)$
  \be\la{1.9}
  \ba  \bar x^a \rho_0\in   L^1 \cap H^1\cap W^{1,q}  ,\quad
\na u_0 \in  L^2 , \quad \rho_0^{1/2}u_0\in L^2, 
   \ea
  \ee with
  \be\la{2.07} \bar x\triangleq(e+|x|^2)^{1/2}\log^{1+\eta_0} (e+|x|^2),\quad\eta_0=\frac{3}{8}- \frac{1}{2\beta}>0.\ee Then  the problem  \eqref{n1}-\eqref{n4} has a unique global strong solution $(\n,u)$ satisfying that for any $0<   T<\infty,$ \be\la{1.10}\begin{cases}
  \rho\in C([0,T];L^1 \cap H^1\cap W^{1,q} ),\\  \bar x^a\rho\in L^\infty( 0,T ;L^1\cap H^1\cap W^{1,q} ),\\ \sqrt{\n } u,\,\na u,\, \bar x^{-1}u,
   \,    \sqrt{t} \sqrt{\n}  u_t \in L^\infty(0,T; L^2 ) , \\ \na u\in  L^2(0,T;H^1)\cap  L^{(q+1)/q}(0,T; W^{1,q}), \\ \sqrt{t}\na u\in L^2(0,T; W^{1,q} )   ,  \\ \sqrt{\n} u_t, \, \sqrt{t}\na u_t ,\,  \sqrt{t} \bar x^{-1}u_t\in L^2(\rr\times(0,T)) ,
   \end{cases}\ee   and that \be \la{l1.2}\inf\limits_{0\le t\le T}\int_{B_{N }}\n(x,t) dx\ge \frac14\int_{\rr} \n_0(x)dx,\ee
for some constant $N >0$ and $B_{N }\triangleq\left.\left\{x\in\rr\right|
\,|x|<N \right\} .$

\end{theorem}

 If the initial data $(\n_0,m_0)$ satisfy some additional regularity and compatibility conditions, the global strong  solutions become classical ones,   that is,
\begin{theorem}\la{t1} Suppose that \eqref{bet} holds. In addition to  \eqref{1.9},  assume   that  $(\n_0,u_0)$ satisfies
  \be\la{1.c1}
  \begin{cases}  \na^2 \n_0,\,\na^2 \lm(\n_0 ),\,\na^2 P(\n_0 )\in L^2\cap  L^q,\\ \bar x^{\de_0}\na^2   \rho_0  ,\, \bar x^{\de_0}\na^2  \lm(\rho_0 ) ,\, \bar x^{\de_0}\na^2 P(  \rho_0 )  \in  L^2    , \quad  \na^2 u_0 \in  L^2 ,
   \end{cases}\ee   for some     constant $\de_0\in (0,1), $  and the following  compatibility condition:
\be \la{co2}- \mu\lap u_0 - \nabla((\mu + \lam(\n_0))\div u_0)+  \nabla P(\n_0)=\n_0^{1/2}g , \ee   with  some $g\in L^2 .$  Then,   in addition to \eqref{1.10} and  \eqref{l1.2},    the  strong  solution  $(\rho,u)$ obtained by Theorem \ref{t1} satisfies  for
  any $0<  T<\infty,$  \be\la{1.a10}\begin{cases}
  \na^2\rho, \,\, \na^2\lm(\rho), \,\,\na^2 P(\rho)\in C([0,T];L^2\cap L^q  ), \\ \bar x^{\de_0}\na^2  \rho   ,\,\, \bar x^{\de_0} \na^2  \lm(\rho )   ,\,\, \bar x^{\de_0} \na^2  P(  \rho)     \in L^\infty( 0,T ;L^2 ) ,\\     \na^2 u,\, \sqrt{\n}  u_t, \,\sqrt{t}   \na  u_t,\,\sqrt{t}  \bar x^{-1}  u_t,\, t\sqrt{\n}u_{tt}, \,t  \na^2 u_t\in L^\infty(0,T; L^2),\\ t\na^3 u\in  L^\infty(0,T; L^2\cap L^q), \,\\    \na u_t,\, \bar x^{-1}u_t,\,  t\na u_{tt},\, t\bar x^{-1}u_{tt}\in L^2(0,T;L^2), \\   t \na^2(\n u)\in L^\infty(0,T;L^{(q+2)/2}) .
   \end{cases}\ee    \end{theorem}

A few remarks are in order:

\begin{remark}  As shown by \cite[Remark 1.1]{hlma},  the   solution $(\n,u )$ obtained in Theorem \ref{t1} is  in fact a classical
    one to  the Cauchy problem  (\ref{n1})-(\ref{n4})
     in $\O\times (0,\infty).$
 \end{remark}

\begin{remark}   Theorems \ref{t2} and \ref{t1} generalize  and improve  the earlier results due to Vaigant-Kazhikhov \cite{Ka} where they required that  $\beta>3$  and that the domain is bounded.  Moreover,  Theorems \ref{t2} and \ref{t1} also extend our previous result   \cite{hlia} where we consider the periodic case to the Cauchy problem in the whole space $\O$. \end{remark}

\begin{remark} It is worth noting here that Zhang-Fang   \cite[Theorem 1.8]{zti} showed that if $(\n, u) \in C^1([0, T ];H^k), k > 3 $ is a spherically
symmetric solution to the Cauchy problem \eqref{n1}-\eqref{n4}  with the compact supported  initial
density $\n_0\not\equiv 0,$      $T$ must be finite provided $1<\beta\le\ga.$  However, in our Theorem \ref{t1},   for $\n,$  we have $\n\in C([0,T];H^2),$  but for $u,$   only  $\na u\in H^k.$ Note that the function $u\in \{\na u\in H^k\}$ decays much slower   for large values of the spatial variable $x$ than the one $u\in H^{k+1}.$ Therefore, it seems that it is the slow  decay of the velocity field for large values of the spatial variable $x$   that leads to the global existence of smooth solutions.      \end{remark}

\begin{remark} It should be mentioned here that it seems that  $\beta>1$  is the extremal case for the system \eqref{n1}-\eqref{n3} (see \cite{Ka}). Therefore, it would be interesting to study the problem \eqref{n1}-\eqref{n4}  when $1< \beta\le 4/3.$ This is left for the future.  \end{remark}

We now comment on the analysis of this paper. Note that for  initial
data satisfying the conditions of Theorems \ref{t2} and \ref{t1},  the local existence and uniqueness of  strong and classical
solutions to the Cauchy problem  (\ref{n1})-(\ref{n4})  have been
established  in \cite{hlma}. Thus, to extend the strong and classical
solutions globally in time, one needs global a priori estimates   on
smooth solutions to (\ref{n1})-(\ref{n4}) in suitable higher norms. To do so,  motivated by   \cite{hlx,hlia},  it turns out
that the key issue in this paper is to derive  the
  upper bound for the density.    We then try to modify the analysis in \cite{Ka,hlia}. However, the methods in  \cite{Ka,hlia} can not be applied directly to our case since  their arguments   rely  heavily on the fact that  the domain is  bounded.  
The key steps of this paper are as follows: We first obtain  the spatial weighted mean estimate  of the density(see \eqref{o3.7}).  Then,  rewriting $\eqref{n1}_2$ as \eqref{key}
 in terms of  a sum of  commutators of Riesz transforms and the operators of multiplication by $u_i$ (see \eqref{a3.42}) as in \cite{L2,Mik,hlia},     we succeed in deriving the estimate of $L^\infty(0,T;L^p)$-norm of the density (see \eqref{ja1})   after using  the spatial weighted mean estimate  of the density we  have just derived, the Hardy type inequality (see \eqref{lpg}),   and the $L^p$-estimate of the commutators due to  Coifman-Rochberg-Weiss \cite{coi1}(see \eqref{2.6}). Next, by energy type estimates   and the compensated  compactness analysis
  \cite[Theorem II.1]{coi3},  we  show that $\log (1+\|\na u\|_{L^2})$
 does not exceed   $\|\n\|^{4/3}_{L^\infty}$ (see \eqref{n9}).
  Then, after  we establish  a key estimate of $\|\n u\|_{L^r}$  in terms of $\|\n\|_{L^\infty},$ $\|\n^{1/2} u\|_{L^2} ,$ and   $\|\na u\|_{L^2} $  with the explicit expression    of $ r $ (see \eqref{2.16} for details),    we can
 use the $W^{1,p}$-estimate  of the   commutator  due to  Coifman-Meyer \cite{coi2} (see \eqref{2.7})  to obtain an estimate on the $L^1(0,T;L^\infty)$-norm of the commutators in terms of $\|\n\|_{L^\infty}$ (see \eqref{a3.41}), which  together with  the Brezis-Wainger inequality   (see \eqref{bmo})  leads to the key a priori
estimate on $\|\n\|_{L^\infty}$    provided $\beta>4/3.$  See Proposition \ref{aupper} and its proof.

 The
next main step is to bound the gradients of the density.   We first obtain the temporal weighted mean estimates on the material derivatives of the velocity  by modifying the basic estimates on the material derivatives of the velocity  due to Hoff \cite{Ho4}. Then, following  \cite{hlx},   the $L^p$-bound of the gradient  of the density can be obtained by solving a logarithm
Gronwall inequality based on a Beale-Kato-Majda type inequality
(see Lemma \ref{le9}),   the a priori estimates we have just
derived and some careful  initial layer analysis;   and  moreover, such a derivation yields simultaneously
also the bound for $L^1(0,T;L^\infty({\O} ))$-norm of the
gradient of the velocity; see Lemma \ref{le5} and its proof.

The rest of the paper is organized as follows: In Section 2, we collect some
elementary facts and inequalities which will be needed in later analysis. Section 3
is devoted to the derivation of  upper    bound  on the density   which  is  the key to extend the local solution to all time. Based on the previous estimates, higher-order ones  are established in Sections 4 and 5. Then finally, the main results,
Theorems \ref{t2} and \ref{t1}, are proved in Section 6.

\section{Preliminaries}

The following    local existence of strong and classical solutions  can be found in \cite{hlma}.

\begin{lemma}   \la{th0}  Let $\beta\ge 1$ and  $\ga>1.$    Assume  that
 $(\n_0,u_0 )$ satisfies \eqref{1.9}.   Then there exist  a small time
$T >0$    and a unique strong solution $(\rho , u )$ to the
  problem   \eqref{n1}-\eqref{n4}  in
$\O\times(0,T )$ satisfying \eqref{1.10} and  \eqref{l1.2}. Moreover, if  $(\n_0,u_0)$ satisfies \eqref{1.c1} and    \eqref{co2} besides \eqref{1.9}, $(\rho , u )$  satisfies \eqref{1.a10} also.
 \end{lemma}

The   following  Sobolev  inequality will be used frequently.
\begin{lemma}[\cite{la,tal}] \la{leo}
    There exists a universal
    positive constant $C$
     such that the following estimates hold  for any $p\in (2,\infty),$
       \be\la{lp}       \|u\|_{L^p }\le C  p^{1/2}      \|\na u\|_{L^{2p/(p+2)}}, \quad \|v\|_{L^p }\le C  p^{1/2} \|v\|_{L^2 }^{2/p}     \|\na v\|_{L^2}^{1-2/p}   , \ee for any function $u\in \left\{ u\in L^p\left|\na u\in L^{2p/(p+2)}\right.\right\}$  and $v\in H^1.$
\end{lemma}

The following weighted $L^p$ bounds for elements of the Hilbert space $\ti D^{1,2}(\O)\triangleq\{u\in H^1_{\rm loc}(\O)|\na u\in L^2(\O)\} $ can be found in \cite[Theorem B.1]{L2}.
\begin{lemma} \la{1leo}
   For   $m\in [2,\infty)$ and $\theta\in (1+m/2,\infty),$ there exists a positive constant $C$ such that we have for all $v\in \ti D^{1,2}(\O),$ \be\la{3h} \left(\int \frac{|v|^m}{e+|x|^2}(\log (e+|x|^2))^{-\theta}dx  \right)^{1/m}\le C\|v\|_{L^2(B_1)}+C\|\na v\|_{L^2 },\ee  where we denote $B_N\triangleq\{x\in\O|
\,|x|<N\} $ for $N>0.$
\end{lemma}

The combination of Lemma \ref{leo} with  Lemma \ref{1leo}   yields
\begin{lemma}  \la{leg}
  For $\bar x$ and $\eta_0$ as in \eqref{2.07},  there exists a
    positive constant $C $ depending only on $\eta_0$
     such that    every function $v\in \ti D^{1,2}(\O)$   satisfies for all $ \de\in (0,2)$\be \la{lpg} \|   v\bar x^{-\de} \|_{L^{4/\de}} \le  C  \de^{-\eta_0-1/2}  \left( \|v\|_{L^2(B_1)}+ \|\na v\|_{L^2} \right).\ee
\end{lemma}

{\it Proof.} Noticing that
\bnn |\pa_i\bar x|\le 6  \log^{1+\eta_0} (e+|x|^2) ,\quad i=1,2,\enn
 we obtain by direct calculations
\be\la{2.h10}\ba &\|\na (v\bar x^{-\de})\|_{L^{4/(2+\de)}}\\ &=\|\bar   x^{-\de}\na v -\de  v\bar x^{-\de-1}\na\bar x \|_{L^{4/(2+\de)}}\\ &\le  \|\na v\|_{L^2}\|\bar   x^{-\de}\|_{L^{4/\de}}+6\de\| v\bar x^{-1}\|_{L^2}\|\bar   x^{-\de}\log^{1+\eta_0}(e+|x|^2)\|_{L^{4/\de}}\\&\le C \left(\|\na v\|_{L^2}+  \de^{-\eta_0}\| v\bar x^{-1}\|_{L^2}\right),\ea\ee
where in the last inequality we have used the following simple fact that
\bnn\ba &\|\bar   x^{-\de}\log^{1+\eta_0}(e+|x|^2)\|_{L^{4/\de}} \\ &\le \|\bar x^{-\de/2}\|_{L^{4/\de} }\|(e+|x|^2)^{- \de/(4(1+\eta_0))}\log(e+|x|^2)\|_{L^\infty}^{1+\eta_0}\\ &\le C \de^{-1-\eta_0} ,\ea\enn due to   $ (e+y)^{-\alpha} \log (e+y)\le \alpha^{-1}$ for $\alpha>0$ and any $y\ge 0.$ The desired estimate \eqref{lpg} thus directly follows from  \eqref{lp}, \eqref{2.h10}, and \eqref{3h}. The proof of Lemma \ref{leg} is completed.

A useful consequence of Lemma \ref{leg} is the following weighted  bounds for elements of  $\ti D^{1,2}(\O) $ which   is important for our analysis.
\begin{lemma}\la{lemma2.6} Let  $\bar x$ and $\eta_0$ be as in \eqref{2.07}. For  $\ga>1,$ assume that $\n \in L^1(\O)\cap L^\infty(\O)$ is a non-negative function such that
\be \la{2.12}   \int_{B_{N_1} }\n dx\ge M_1,  \quad \int \n^\ga dx\le M_2, \quad  \int \n \bar x^\alpha dx\le M_3,\ee
for positive constants $   M_i (i=1,\cdots,3),\alpha ,$ and $ N_1\ge 1.$  Then there is a positive constant $C$ depending only on   $ M_i (i=1,\cdots,3), N_1,\alpha, \ga,$ and $\eta_0$ such that      every $v\in \ti D^{1,2}(\O)  $ satisfies
 \be\la{2.16}\|\n v\|_{L^r}   \le Cr^{\eta_0+1/2}(1+ \|\n\|_{L^\infty})   \left(  \|\n^{1/2} v\|_{L^2(B_{N_1})}
 + \|\na  v\|_{L^2}\right) ,\ee
  for any $r\in (1,\infty).$

\end{lemma}

{\it Proof.}   It follows from   \eqref{2.12} and  the Poincar\'e type inequality     \cite[Lemma 3.2]{Fe}  that  there exists a positive constant $C $ depending only on $  M_1, M_2, N_1 ,$ and $\ga,$  such that  \be \la{3.12} \|v\|_{H^1(B_{ N_1} )}^2\le C \int_{B_{ N_1} }\n v^2dx +C \|\na v\|_{L^2(B_{ N_1} )}^2.\ee  This combined with  Holder's inequality,  \eqref{lpg}, and   \eqref{2.12} yields that   for  $r\in (1,\infty) $  and $\si=4/(r(4+\alpha)),$
 \bnn\ba \|\n v\|_{L^r}  &\le   \|(\n \bar x^\alpha)^\si\|_{L^{1/\si}} \| v\bar x^{-\alpha\si}\|_{L^{4/(\alpha\si)}} \|\n\|_{L^\infty}^{1-\si}\\ &\le Cr^{\eta_0+1/2}   \left(  \|\n^{1/2} v\|_{L^2(B_{N_1})}
 + \|\na  v\|_{L^2}\right)(1+ \|\n\|_{L^\infty}) ,\ea\enn
 which shows \eqref{2.16} and  finishes the proof of Lemma \ref{lemma2.6}.

Next,   let $\mathcal{H}^1(\O)$ and $\mathcal{BMO}(\O)$ stand for the usual {\sc Hardy} and {\sc BMO} space.
Given a function $b,$ define the linear operator
$$[b,R_iR_j](f)\triangleq bR_i\circ R_j(f)-R_i\circ R_j (bf), \,\,i,j= 1, 2,$$ where   $R_i$ is the usual
 {\sc Riesz} transform on $\O: R_i = (-\lap)^{-1/2}\p_i$.
The following properties of the  commutator $[b,R_iR_j](f)$   will be   useful for our analysis. \begin{lemma} Let $b,f\in C_0^\infty(\O).$  Then for $p\in (1,\infty),$ there is   $ C(p)$ such that
\be\la{2.6} \|[b,R_iR_j](f)\|_{L^p}\le C(p)\|b\|_{\mathcal{BMO}}\|f\|_{L^p}.
\ee
Moreover, for $q_i\in (1,\infty) (i=1,2,3)$ with $q_1^{-1}=q_2^{-1} +q_3^{-1},$ there is  x $ C$ depending only on $q_i (i=1,2,3)$ such that
\be \la{2.7}\|\na[b,R_iR_j](f)\|_{L^{q_1}}\le C \|\na b\|_{L^{q_2}}\|f\|_{L^{q_3}}.
\ee
\end{lemma}

\begin{remark} Properties \eqref{2.6} and \eqref{2.7} are  due to Coifman-Rochberg-Weiss \cite{coi1} and Coifman-Meyer \cite{coi2} respectively.
\end{remark}

Next, we state the following Beale-Kato-Majda-type inequality
which was proved in \cite{B1,kato} when $\div u\equiv 0$ and will be used
later to estimate $\|\nabla u\|_{L^\infty}$ and
$\|\nabla\rho\|_{L^p}$.
\begin{lemma}[\cite{B1,kato}]   \la{le9}  For $2<q<\infty,$ there is a
constant  $C(q)$ such that  the following estimate holds for all
$\na u\in  L^2(\O)\cap W^{1,q} (\O),$ \bnn\la{ww7}\ba \|\na u\|_{L^\infty
}&\le C\left(\|{\rm div}u\|_{L^\infty }+ \|{\rm rot} u\|_{L^\infty }
\right)\log(e+\|\na^2 u\|_{L^q })+C\|\na u\|_{L^2 } +C .
\ea\enn
\end{lemma}

Finally, the   following     Brezis-Wainger  inequality will also be used.
\begin{lemma}[\cite{en,eng}] \la{gleo} For   $q>2,$     there exists some positive constant $C$ depending only on $q$   such that  every function $v\in \left. \left\{v\in  W^{1,q} \right|\na v\in L^2\right\} $ satisfies
\be\la{bmo}
\norm[L^\infty]{v}\le C  (\|v\|_{L^q}+\|\na v\|_{L^2} )\ln^{1/2}(e+\|  v\|_{W^{1,q}}) +C.
\ee
\end{lemma}
\section{\la{se3}A priori estimates (I): upper bound of the density}

In this section and the next, in addition to the conditions of Theorem \ref{t2}, we will always assume that smooth $(\n_0,u_0)$  satisfies    \be\la{ta2} \n_0(x)>0,  \quad \frac12\le \int_{B_{N_0}}\n_0(x)dx \le \int_{\rr }\n_0(x)dx \le 2,\ee
for some a positive constant $N_0.$
Moreover, suppose that $(\rho,u)$ is the strong solution to  \eqref{n1}-\eqref{n4} on $\O\times (0,T]$ obtained by Lemma \ref{th0}.

The following Proposition \ref{aupper}  will give an  upper bound of the density which is the key to obtain higher order estimates.

\begin{proposition}\la{aupper} Under the conditions of Theorem \ref{t2},  for \bnn E_0\triangleq  \|\n_0\bar x^\al\|_{L^1}+ \|\n_0\|_{L^\infty}+ \|\n_0^{1/2}u_0 \|_{L^2}+ \|\na u_0\|_{L^2}, \enn there is a positive constant  $C $ depending only on $  \mu,  \beta,$ $ \gamma,   T, N_0, \al,$ and $E_0$ such that
   \be\la{b3.56}
\sup_{0\le t\le T}\left(\|\n\|_{L^\infty}+\|\na u\|_{L^2}\right)  + \int_0^T\int \n | u_t+u\cdot\na u|^2dxdt\le C .
\ee
\end{proposition}

 Before proving Proposition \ref{aupper}, we establish some a priori estimates, Lemmas \ref{kq2}--\ref{lemma1}.

First, we have the following lemma.

\begin{lemma}  \la{kq2}
There exist positive constants $C $ and $N_1$ both depending only on $\al, \ga, T, N_0,  $   $ \|\n_0\bar x^\al\|_{L^1},   \|\n_0\|_{L^\ga},$ and $ \|\n_0^{1/2}u_0 \|_{L^2}$ such that
\be\la{r3}
 \sup\limits_{0\le t\le T}\int\left(\n|u|^2+\n^\ga+\n\bar x^\al \right)dx +\int_0^T \int\left( \mu |\na u|^2+ \lambda(\n) (\div u)^2\right)dxdt \le C,
\ee and \be\la{rr4}\ba \inf\limits_{0\le t\le T}\int_{B_{N_1} }\n  dx \ge 1/4.\ea\ee
\end{lemma}

{ \it Proof. }  First,
  the standard energy inequality reads:
\be\la{rr3}
 \sup\limits_{0\le t\le T}\int\left(\n|u|^2+\n^\ga\right)dx +\int_0^T \int\left( \mu |\na u|^2+(\mu+\lambda(\n))(\div u)^2\right)dxdt \le \ti C.
\ee

Next, multiplying $\eqref{n1}_1$ by $\bar x^\al $ and integrating the resulting equality over $\O,$ we obtain after integration by parts and using \eqref{rr3} that
\bnn \ba  \frac{d}{dt}\int\n \bar x^\al  dx &\le C\int \n |u|\bar x^{\al -1}\log^{1+\eta_0}(e+|x|^2)dx\\ &\le C\left(\int\n \bar x^{2\al -2} \log^{2(1+\eta_0)}(e+|x|^2)dx\right)^{1/2}\left(\int\n u^2 dx\right)^{1/2}\\ &\le C\left(\int\n \bar x^{ \al  } dx\right)^{1/2},\ea\enn
which together with Gronwall's inequality gives
\be\la{o3.7} \sup_{0\le t\le T}\int\n \bar x^\al  dx\le C .\ee
This, along with \eqref{rr3}, gives  \eqref{r3}.

Finally, the mass conservation equation  $\eqref{n1}_1$ yields \be \la{mr3}\int\n dx =\int \n_0dx .\ee For $N>1,$ let $\varphi_N$   be a smooth function such that  \be\la{vl1}0\le\varphi_N(x)\le 1,\quad \varphi_N=\begin{cases} 1& \mbox{ if }\,|x|\le N,\\ 0& \mbox{ if }\,|x|\ge 2N,\end{cases} \quad |\na \varphi_N|\le 2N^{-1}.\ee
 It follows from $(\ref{n1})_1,$ \eqref{mr3}, and \eqref{rr3} that  \bnn\ba \frac{d}{dt}\int \n \vp_N dx &=\int \n u \cdot\na \vp_N dx\\ &\ge - 2N^{-1}\left(\int\n dx\right)^{1/2}\left(\int\n |u|^2dx\right)^{1/2}\ge -2\ti C^{1/2}N^{-1},\ea\enn which gives   \bnn  \inf\limits_{0\le t\le T}\int \n \vp_N dx\ge \int \n_0 \vp_N dx-2\ti C^{1/2}N^{-1}T.\enn
  This combined with   \eqref{ta2} yields that for  $N_1\triangleq 2(2+ N_0+8\ti C^{1/2}T ),$
\bnn\ba  \inf\limits_{0\le t\le T}\int_{B_{N_1}  }\n  dx\ge \int \n \vp_{N_1/2}dx\ge 1/4,\ea\enn
 which shows \eqref{rr4}.   The proof of Lemma \ref{kq2} is completed.

\begin{lemma} \la{kq1}  Assume that \eqref{bet} holds. Then there is a positive constant  $C $   depending only on $  \mu,  \beta,$ $ \gamma,   T, N_0, \al,$ and $E_0$   such that
\be\la{ja1}\ba   \sup_{0\le t\le T}\int \left(\n+\n^{2\beta\ga+1 }\right) dx  \le C.\ea\ee
  \end{lemma}

{\it Proof.} First,
  we denote \bnn \nabla^{\perp}\triangleq (\p_2,-\p_1),\quad \frac{D}{Dt}f \triangleq\dot{f}\triangleq f_t  + u \cdot\nabla f ,\enn where $\dot f$ is the material derivative of $f.$ Let $G$ and $\o$ denote the effective viscous flux and the vorticity respectively as follows:
\bnn
G \triangleq (2\mu+\lam(\rho))\div u -   P ,\quad \o \triangleq \nabla^{\perp}\cdot u= \p_2u_1 - \p_1u_2.
\enn
We thus rewrite
the momentum equations $(\ref{n1})_2$ as
\be\la{n5}
\rho\dot{u} = \nabla G + \mu\nabla^{\perp}\o,
\ee
which shows that  $G$ solves
\bnn
\lap G = \div(\rho\dot{u})=\p_t(\div(\rho u)) + \div\div(\rho u\otimes u).
\enn
This implies
\be\la{n12}
\ba
  &G  + \frac{D}{Dt}\left((-\lap)^{-1}\div(\rho u)\right)=F,
\ea
\ee with  the  commutator $F$  defined   by \be \la{a3.42}F\triangleq \sum_{i,j=1}^2[u_i,R_iR_j](\rho u_j)= \sum_{i,j=1}^2\left(u_iR_i\circ R_j(\rho u_j)-R_i\circ R_j (\rho u_iu_j)\right).\ee

  Then, since $\n>0$ due to \eqref{ta2},  the mass equation $(\ref{n1})_1$ leads to \bnn   -\div u =\frac{1}{\n}D_t\n ,\enn which
  combined with (\ref{n12})  gives that
\be\la{key}
 \frac{D}{Dt}\theta(\rho)  + P = \frac{D}{Dt}\psi -F ,
\ee
with \be\la{p26}
\theta(\n) \triangleq 2\mu\log\rho +  \beta^{-1}\rho^{\beta}  ,\quad
 \psi\triangleq(-\lap)^{-1} \div(\rho u).\ee

Next,
 denoting $f\triangleq\max\{\theta(\n)-\psi,0\},$   multiplying  \eqref{key} by $\n f^{2\ga-1},$ and  integrating the resulting equality over $\O$ lead  to
  \be\la{2.11}\ba \frac{d}{dt}\int \n f^{2\ga}dx&\le C\int \n f^{2\ga-1}|F|dx\\ & \le C\|\n^{1/(2\ga)}f\|^{2\ga-1}_{L^{2\ga}} \|\n\|_{L^{2\beta\ga+1}}^{1/(2\ga)}
  \|F\|_{L^{(2\beta\ga+1)/\beta}}\\ & \le C\|\n^{1/(2\ga)}f\|^{2\ga-1}_{L^{2\ga}} \|\n\|_{L^{2\beta\ga+1}}^{1/(2\ga)}
  \|\na u\|_{L^2}\|\n u\|_{L^{(2\beta\ga+1)/\beta}} ,
  \ea\ee where in the last inequality we have used the following simple fact that for any $p\in (1,\infty),$
 \be\la{2.9}\|F\|_{L^p}\le C(p)\|u\|_{\mathcal{BMO}}\|\n u\|_{L^p}\le C(p)\|\na u\|_{L^2}\|\n u\|_{L^p},\ee due to \eqref{2.6}.
      It follows from the Holder inequality, \eqref{r3}, \eqref{lpg},  \eqref{rr4}, and \eqref{3.12}   that
\be\la{2.10}\ba  \|\n u\|_{L^{(2\beta\ga+1)/\beta}} &\le \|(\n\bar x^\al)^\si\|_{L^{1/\si }}\|u \bar x^{-\al\si} \|_{L^{4/(\al\si) }}\|\n^{1-\si}\|_{L^{(2\beta\ga+1)/(1-\si)}}\\ &\le C  (1+\|\na u   \|_{L^2}) \left(1+\|\n\|_{L^{2\beta\ga+1}}\right),\ea\ee
 where
$\si=4(\beta-1)/((4+\al)(2\beta\ga+1)-4).$
   Substituting \eqref{2.10} into \eqref{2.11} gives
  \be\la{3.25}\ba \frac{d}{dt}\int \n f^{2\ga}dx   \le C\left(1+\int \n f^{2\ga}dx+ \int \n^{2\beta\ga+1}  dx  \right)(1+
  \|\na u\|_{L^2}^2),
  \ea\ee
 due to $\beta>1.$

Next, we claim
  \be\la{3.29}\ba  \int \n^{2\ga\beta+1}dx \le C+C\int \n f^{2\ga }dx, \ea\ee
 which together with \eqref{3.25}, \eqref{r3},  and Gronwall's inequality yields \bnn  \sup_{0\le t\le T}\int \n f^{2\ga}dx\le C.\enn
This combined with \eqref{3.29} and \eqref{mr3} directly gives  \eqref{ja1}.

 Finally, it only remains to prove \eqref{3.29}. In fact, for any $p\in (1,\infty),$ we have \be \la{3.a29} \|\na \psi\|_{L^p}\le C(p)\|\n u\|_{L^p},\ee which, along with \eqref{lp} and \eqref{r3}, gives that for any $r\in (2,\infty),$ \be \la{f3.18} \ba\|\psi\|_{L^r}&\le C(r)\|\na\psi\|_{L^{2r/(r+2)}}\\&\le C(r)
    \|\n^{1/2} \|_{L^{r}}\|\n^{1/2}u\|_{L^2}\le C(r)
    \|\n\|_{L^{r/2}}^{1/2}.\ea \ee It thus follows from  \eqref{f3.18}, \eqref{mr3}, and  the fact that $\beta>1  $ that
  \bnn\ba \int \n^{2\ga\beta+1}dx  &= \int_{(\n\le 2)}\n^{2\ga\beta+1}dx+\int_{(\n>2)}\n^{2\ga\beta+1}dx\\ &\le C\int_{(\n\le 2)} \n dx+C\int_{(\n>2)}\n f^{2\ga }dx+C\int \n |\psi|^{2\ga}dx \\ &\le C +C\int \n f^{2\ga }dx+C\|\n\|_{L^{(2\beta\ga+1)/(2\beta\ga+1-\ga)}} \|\psi\|_{L^{2(2\beta\ga+1)}}^{2\ga}\\ &\le C +C\int \n f^{2\ga }dx+C(1+\|\n\|_{L^{ 2\beta\ga+ 1 }} ) \|\n \|_{L^{ 2\beta\ga+1 }}^{ \ga}\\ &\le C(\ve) +C\int \n f^{2\ga }dx+ \ve\int \n^{2\ga\beta+1}dx .\ea\enn
     The proof of Lemma \ref{kq1} is finished.

 The following  $L^p$-estimate of the momentum will play  an important role in the estimate of the upper bound of the density.

\begin{lemma}  \la{jle2} Assume that \eqref{bet} holds. Then,  for any $p>4,$  there is a positive constant  $C(p)$ depending only on $ p, \mu,  \beta,$ $ \gamma,   T, N_0, \al,$ and $E_0$ such that
\be\la{jle3}\ba \|\n u\|_{L^p}   \le  C(p) \rt^{1+\beta/4+(\beta\eta_0)/2}  (1+ \|\na u\|_{L^2})^{ 1-2/p}, \ea\ee with $\eta_0$    as in \eqref{2.07} and \bnn  \rt \triangleq 1+\sup\limits_{0\le t\le T} \|\n\|_{L^\infty}.\enn
  \end{lemma}

 {\it Proof. } First, for $$\nu\triangleq \frac{\mu^{1/2}}{ 2(\mu+1) } \rt^{-\beta/2}\in (0,1/4],$$  multiplying $(\ref{n1})_2$ by $(2+\nu)|u|^\nu u$  and   integrating the resulting equation over $ \O$ lead  to
\bnn\ba & \frac{d}{dt}\int \n |u|^{2+\nu}dx+ (2+\nu) \int|u|^\nu \left(\mu |\na u|^2+(\mu+\n^\beta) (\div u)^2\right) dx \\& \le  (2+\nu)\nu  \int (\mu+\n^\beta)|\div u||u|^\nu |\na u|dx +C  \int \n^\ga |u|^\nu |\na u|dx\\& \le \frac{2+\nu}{2}\int (\mu+\n^\beta)(\div u)^2|u|^\nu  dx +\frac{(2+\nu)\mu}{8 (\mu+1)} \int  |u|^\nu |\na u|^2 dx \\&\quad+ \mu  \int   |u|^\nu |\na u|^2dx + C \int \n |u|^{2+\nu} dx+C\int \n^{(2+\nu)\ga-\nu/2}dx ,\ea\enn
which  together with Gronwall's inequality and \eqref{ja1} thus gives
   \be\la{ja2}
  \sup_{0\le t\le T}\int \n |u|^{2+\nu}dx\le C .
  \ee

Then, it follows from Holder's inequality,   \eqref{ja2},   \eqref{r3},   \eqref{rr4},  and  \eqref{2.16}    that for $ r=(p-2)(2+\nu)/\nu,$
\bnn\ba \|\n u\|_{L^p}  &\le  \|\n u\|_{L^{2+\nu}}^{2/p}\|\n u\|_{L^r}^{1-2/p
}\\ &\le C \rt^{ (1+\nu)/p}  \|\n^{1/(2+\nu)} u\|_{L^{2+\nu}}^{2/p} \left(r^{\eta_0+1/2}\rt (1+ \|\na  u\|_{L^2})\right)^{ 1-2/p}\\ &\le C(p)\rt^{ (1+\nu)/p}  \left(\rt^{1+\beta/4+(\beta\eta_0)/2}(1+ \|\na  u\|_{L^2})\right)^{ 1-2/p} \\ &\le  C(p)\rt^{1+\beta/4+(\beta\eta_0)/2} (1+ \|\na  u\|_{L^2})^{ 1-2/p},\ea\enn  which   shows \eqref{jle3} and finishes   the proof of  Lemma \ref{jle2}.

\begin{lemma}\la{le1}Assume that \eqref{bet} holds.
   Then   there is a constant $C   $ depending only on $  \mu,  \beta,$ $ \gamma,   T, N_0, \al,$ and $E_0$  such that
\be\la{n9}
\ba
  &\sup\limits_{0\le t\le T}\log(e+A^2(t)) + \int_0^T\frac{B^2(t)}{e + A^2(t)}dt    \le   C \rt^{4/3},
\ea
\ee
where
\be\la{an7}
 A^2(t) \triangleq \int \left( \o^2 (t) + \frac{ G^2(t) }{2\mu+\lam(\n(t))}\right)dx,\quad  B^2(t) \triangleq \int \rho(t)|\dot{u}(t)|^2dx.\ee

\end{lemma}

{\it Proof.}  First,  direct calculations show that
\be\la{ra3}\ba
  \nabla^{\perp}\cdot \dot u = \frac{D}{Dt}\o -(\p_1u\cdot\na) u_2+(\p_2u\cdot\na)u_1  = \frac{D}{Dt}\o + \o \div u , \ea
\ee
and that
\be\la{rb3}\ba \div  \dot u&=\frac{D}{Dt}\div u +(\p_1u\cdot\na) u_1+(\p_2u\cdot\na)u_2 \\&=
\frac{D}{Dt} \frac{G}{2\mu+\lam} + \frac{D}{Dt} \frac{P}{2\mu+\lam}  - 2\nabla u_1\cdot\nabla^{\perp}u_2 + (\div u)^2.\ea
\ee

Then, multiplying  \eqref{n5}    by $2\dot u $  and integrating the resulting equality over $\O,$  we obtain after using \eqref{ra3} and \eqref{rb3} that
\be\la{p3}
\ba
  \frac{\rm d}{{\rm d}t}A^2  + 2B^2
& = -\int \o^2\div udx + 4\int G\nabla u_1\cdot\nabla^{\perp}u_2dx- 2\int G(\div u)^2dx\\
 &\quad -\int\frac{(\beta-1)\lam - 2\mu}{(2\mu+\lam)^2}G^2\div udx + 2\beta \int\frac{\lam(\n) P}{(2\mu+\lam)^2}G\div udx\\
 & \quad-2\g\int\frac{ P}{2\mu+\lam}G\div udx \triangleq \sum_{i=1}^6I_i.
 \ea
\ee

Each $I_i$ can be estimated as follows:

First,
it follows from \eqref{n5} that \bnn \lap G=\div(\n \dot u),\quad \mu\lap \o=\na^\perp\cdot (\n\dot u),\enn
which together with the standard $L^p$-estimate of elliptic equations yields that for $p\in (1,\infty),$\be \la{3.36}\|\na G\|_{L^p}+\|\na \o\|_{L^p}\le C(p,\mu)\|\n \dot u\|_{L^p}.\ee  In particular, we have
\be\la{d3.19}  \|\na G\|_{L^2}+\|\na \o\|_{L^2}\le C(\mu)\rt^{1/2}B.\ee
This combined with \eqref{lp} gives
\be\la{a3.26}
\ba   \|\o\|_{L^4} & \le  C\|\o\|_{L^2}^{1/2}\|\na\o\|_{L^2}^{1/2}\\&\le   C \rt^{1/4}
A^{1/2} B^{1/2},
\ea\ee
which leads to
\be\la{a3.27}\ba
|I_1|  \le C\|\o\|_{L^4}^2\|\div u\|_{L^2}     \le \ve B^2+C(\ve)\rt \|\na u\|_{L^2}^2 A^2 .\ea
\ee

Next, we will use an idea due to  \cite{des,Mik} to estimate $I_2.$ Noticing that $${\rm rot} \na u_1=0,\quad \div \nabla^{\perp}u_2=0,$$ one derives   from \cite[Theorem II.1]{coi3} that $$ \|\na u_1\cdot \nabla^{\perp}u_2\|_{\mathcal{H}^1}\le C\|\na u\|_{L^2}^2.$$
This combined with the fact that $\mathcal{BMO}$ is the dual space of $\mathcal{H}^1$ (see \cite{fef1}) gives
\be \la{a3.28}\ba |I_2|&\le C\|G\|_{\mathcal{BMO}}\|\na u_1\cdot \nabla^{\perp}u_2\|_{\mathcal{H}^1}\\&\le C\|\na G\|_{L^2}\|\na u\|_{L^2}^2  \\&\le C
\rt^{1/2}B  \|\na u\|_{L^2}(1+A) \\&\le \ve B^2+C(\ve)\rt \|\na u\|_{L^2}^2(1+A^2),   \ea\ee
where in the third inequality we have used \eqref{d3.19} and  the following simple fact  that  for $t\in [0,T],$
   \be \la{a3.12}
 C^{-1}\|\na u(\cdot,t)\|_{L^2}^2 -C\le A^2(t)\le  C \rt^{\beta}\|\na u(\cdot,t)\|_{L^2}^2 +C, \ee due to  \eqref{ja1}.

Next, Holder's inequality yields that for $\de\in (0,2(\beta-1)),$
\be\la{pp4}\ba
 \sum\limits_{i=3}^6|I_i|
 &\le C\int|\div u|\left(|G| \frac{|G|+P}{2\mu+\lam}+  \frac{G^2}{2\mu+\lam} +  \frac{P |G| }{2\mu+\lam}  \right)dx \\
 &\le   C\int \frac{G^2 |\div u|}{2\mu+\lam}dx  +C\int \frac{P |G| }{2\mu+\lam} |\div u|dx   \\& \le  C \|\na u\|_{L^2}  \left\|\frac{G^2}{2\mu+\lam}\right\|_{L^2}+  C  \|\na u\|_{L^2} \|  G\|_{L^{2(2+\de)/\de}} \|  P  \|_{L^{2+\de}} \\& \le  C \|\na u\|_{L^2} \left\|\frac{G^2}{2\mu+\lam}\right\|_{L^2}+  C(\de)  \|\na u\|_{L^2} \|  G\|^{\de/(2+\de)}_{L^2} \| \na G\|^{2/(2+\de)}_{L^2}     ,\ea
\ee where in the last inequality we have used \eqref{ja1} and \eqref{lp}.

Then,   noticing that  \eqref{an7} gives
\be\la{ggaa} \ba \|G\|_{L^2}  \le C \rt^{\beta/2} A  ,\ea\ee  one deduces from
the Holder inequality and \eqref{lp}   that for $0<\de<1,$
\be\la{p7}\ba
\left\|\frac{G^2}{\sqrt{2\mu+\lam}}\right\|_{L^2}&\le C\left\|\frac{G}{\sqrt{2\mu+\lam}}\right\|_{L^2}^{1-\de}
\|G\|_{ L^{2(1+\de)/\de}}^{1+\de}\\ &\le C(\de)A^{ 1-\de }\| G\|_{L^2}^{\de}
\|\na G\|_{L^2}   \\ &\le C(\de)\rt^{ (1+\de\beta)/2 }AB
  ,\ea
\ee where in the last inequality we have used \eqref{d3.19}.
Putting  \eqref{p7}, \eqref{ggaa},   and \eqref{d3.19} into \eqref{pp4}  yields
\be \la{p4}\ba
 \sum\limits_{i=3}^6|I_i|   &\le C(\de) \rt^{ (1+\de\beta)/2 }\|\na u\|_{L^2} \left( A B+    A^{\de/(2+\de)}B^{2/(2+\de)} \right)\\  &\le C(\de) \rt^{ (1+\de\beta)/2 }\|\na u\|_{L^2}  ( A B+B + A  )\\  &\le \ve B^2+C(\ve,\de) \rt^{  1+\de\beta  }(1+\|\na u\|^2_{L^2})(1+A^2)   .\ea
\ee

Finally, substituting \eqref{a3.27},  \eqref{a3.28}, and  \eqref{p4}   into  \eqref{p3},  we obtain after choosing $\ve $ suitably small   that for $\de\in (0,\min\{1,2(\beta-1)\})$ \be
 \la{ou1} \ba
   \frac{d}{dt}A^2  +  B^2   &\le C(\de) \rt^{  1+\de\beta  }(1+\|\na u\|^2_{L^2})(1+A^2).
  \ea
  \ee
 Dividing \eqref{ou1} by $e+A^2,$ choosing $\de= 1/(3\beta),$ and using \eqref{r3}, we obtain \eqref{n9} and finish the proof of Lemma \ref{le1}.

  Next, the following lemma   gives an   estimate of the $L^1(0,T;L^\infty)$-norm of the commutator $F$ defined by \eqref{a3.42}.

\begin{lemma} \la{lemma1}  Assume that \eqref{bet} holds. Then   there is a positive constant  $C $ depending only on  $ \mu,  \beta,$ $ \gamma,   T, N_0, \al,$ and $E_0$ such that
  \be\la{a3.41}
  \ba
   \int_0^T\|F\|_{L^\infty}dt\le   C \rt^{1+\beta/4+ \beta\eta_0 }  .
  \ea
  \ee\end{lemma}

{\it Proof.}
First, it follows from
the Gagliardo-Nirenberg  inequality, \eqref{2.9}, and \eqref{2.7} that for $p\in (8,\infty),$
\be\la{a3.46} \ba\|F \|_{L^\infty}&\le  C(p) \|F \|_{L^p}^{ (p-4)/p}\|\na F\|_{L^{4p/(p+4)}}^{4/p}\\&\le   C(p)\left(\|\na u\|_{L^2} \|\n  u\|_{L^p} \right)^{ (p-4)/p}\left(\|\na u\|_{L^4} \|\n  u\|_{L^p} \right)^{4/p}\\ &\le   C(p)\|\na u\|_{L^2} ^{ (p-4)/p}\|\na u\|_{L^4}^{4/p} \|\n  u\|_{L^p}\\& \le C(p)\rt^{1 +\beta/4+(\beta\eta_0)/2}   \left(1+\|\na u\|_{L^2}\right)^{ 2-6/p} \|\na u\|_{L^4}^{4/p}   ,\ea \ee
where in the last inequality we have used   \eqref{jle3}.

Next,
we obtain from  \eqref{a3.26},  \eqref{p7}, \eqref{a3.12},    and   \eqref{r3}   that
\be\la{3.43}
\ba \|\na u\|_{L^4}   &
 \le C(\|\div u\|_4+\|\omega\|_4)\\ &\le  C\left\|\frac{G+P }{
2\mu+\lambda}\right\|_{L^4} +  C \rt^{1/4}
 A^{1/2} B^{1/2}\\ &\le  C\left\|\frac{G^2 }{
\sqrt{2\mu+\lambda}}\right\|_{L^2}^{1/2} +C \left\|\frac{ P }{
2\mu+\lambda}\right\|_{L^4} +  C \rt^{1/4}
 A^{1/2}B^{1/2}\\& \le C  \rt  A^{1/2}B^{1/2}
 +C\rt^{\ga}    \\
& \le C  \rt^{2\beta\ga }  (e+\|\na u\|_{L^2})  \left(1+\frac{ B^2}{e+A^2}\right)^{1/4}    .\ea\ee
Substituting  \eqref{3.43}  into   \eqref{a3.46} yields   that for  $p>8,$
  \bnn \ba  \|F\|_{L^\infty}   & \le C(p)\rt^{1 +\beta/4 +(\beta\eta_0)/2 +8\beta\ga/p }\left(e+\|\na u\|_{L^2} \right)^{ 2-2/p}  \left(1+\frac{ B^2}{e+A^2}\right)^{1/p}   \\ & \le C(p)\rt^{1 + (\beta/4 +(\beta\eta_0)/2 )  p /(p-1) +9\beta\ga/(p-1) } \left(e+\|\na u\|^2_{L^2}\right) + \frac{ B^2}{e+A^2}, \ea\enn which together with \eqref{n9} and \eqref{r3} directly gives \eqref{a3.41} after choosing $p$ suitably large since $1+\beta/4>4/3$ due to $\beta>4/3.$ The proof of Lemma \ref{lemma1} is completed.

Now we are in a position to prove Proposition \ref{aupper}.

{\it Proof of Proposition \ref{aupper}.}
  For $\psi$   as in \eqref{p26},  it follows from  \eqref{f3.18} and    \eqref{r3} that
  $$\|\psi \|_{L^{2\ga}}\le C,$$
which together with \eqref{bmo},     \eqref{jle3},  and \eqref{3.a29} leads to
\be\la{psi}
\ba
\|\psi\|_{L^\infty} & \le C \left( \|\psi\|_{L^{2\ga}}+ \|{\na\psi}\|_{L^2}\right) \log^{1/2}(e+ \|\psi \|_{W^{1,2\ga}}) +C\\ & \le C\left(1+  \| \n u\|_{L^2} \right) \log^{1/2}(e   +\|\n u \|_{L^{2\ga}}) +C\\ & \le C \rt^{1/2}  \log^{1/2}\left(  \rt^{1+\beta/4+(\beta\eta_0)/2}(e+\|\na u \|_{L^2})\right)+C\\ & \le C \rt^{1/2}\log^{1/2}(e+ A^2) +C\rt\\& \le C \rt^{4/3} ,\ea\ee where in the last inequality we have used \eqref{n9}.
One thus derives from \eqref{key}, \eqref{psi}    and \eqref{a3.41}  that  \bnn    \rt^{\beta}\le  C\rt^{   1+\beta/4+ \beta\eta_0  }. \enn
Because of \eqref{bet} and  \eqref{2.07},  this directly gives \bnn
   \sup\limits_{0\le t\le T}\|\n\|_{L^\infty}\le C ,
  \enn which together with \eqref{ou1},   \eqref{a3.12},    \eqref{r3},   and Gronwall's inequality yields  \eqref{b3.56}.  We complete the proof of Proposition \ref{aupper}.

\section{\la{se4} A priori estimates (II): higher order estimates (I)}

\begin{lemma}\la{310}   Assume that \eqref{bet} holds.  Then there is a positive constant  $C $ depending only on  $   \mu,  \beta,  \gamma,   T, N_0, \al,$ and $E_0$  such that
\be\la{b19}
\sup_{0\le t\le T}t\int \rho|\dot{u}|^2dx + \int_0^Tt\|\na\dot{u}\|_{L^2}^2dt\le C.
\ee
\end{lemma}

{\it Proof.} We will adapt an idea due to \cite{Ho4} to prove \eqref{b19}. In fact, operating $ \pa/\pa t+\div
(u\cdot) $ to $ (\ref{n1})_2^j $   yields that
\be\la{a4.6}\ba &(\n   \dot u_j)_t+\div(\n u  \dot u_j)-\mu \lap \dot u_j-\p_j((\mu+\lm)\div \dot u)\\&=\mu  \p_i(-\p_iu\cdot\na u_j+\div u\p_i u_j)-\mu  \div(\p_i u\p_i u_j)\\&\quad-\pa_j\left[ (\mu+\lambda) \pa_i u\cdot\na u_i-( \mu+(1-\beta)\n^\beta )(\div u)^2\right]\\&\quad-\div (\pa_ju(\mu+\lambda) \div u)+(\ga-1)\p_j(P\div u)+\text{div}(P\p_ju)  .\ea\ee

Then, multiplying \eqref{a4.6} by $\dot u,$
  we obtain after integration
by parts  that
\be\la{b17}
 \ba & \frac{1}{2}\frac{d}{dt}\int \n |\dot u|^2dx+\mu\int |\na \dot u|^2dx+ \int (\mu+\lambda)(\div \dot u)^2dx\\
&\le \frac{\mu}{8}\int|\nabla\dot{u}|^2dx + C \|\na u\|_{L^4}^4 +C \|\nabla   u\|^2_{L^2}    \\
 &\le   \frac{\mu}{8}\int|\nabla\dot{u}|^2dx +C\|\n^{1/2}\dot u\|^2_{L^2} +C,
\ea
\ee where in the second inequality we have used \eqref{3.43} and  \eqref{b3.56}.
Multiplying \eqref{b17} by  $t $  and integrating the resulting inequality  over $(0,T),$ we obtain (\ref{b19}) after using  \eqref{b3.56}.
   We thus finish the proof of Lemma \ref{310}.

\begin{lemma}\la{le5} Assume that \eqref{bet} holds and let $q>2$ be as in Theorem \ref{t2}.     Then there is a constant $C $ depending only on  $ \mu,  \beta,  \gamma,   T, N_0, \al, E_0, q, $ and $\|\n_0\|_{H^1\cap W^{1,q}} $ such that
  \be\la{pa1}\ba
  &\sup_{0\le t\le T}\left(\norm[H^1\cap W^{1,q}]{ \rho} +\|\nabla   u\|_{L^2} +  t\|\na^2 u\|^2_{L^2}  \right)\\&+\int_0^T \left(\|\nabla^2  u\|_{L^2}^2+\|\nabla^2 u\|_{L^q}^{(q+1)/q}+t\|\nabla^2 u\|_{L^q}^2 \right)dt\le C .\ea
  \ee
\end{lemma}

{\it Proof.}  Following \cite{hlx}, we will   prove \eqref{pa1}.
First, denoting    $\Phi\triangleq    (2\mu+\lambda(\n))\na \n ,$ one deduces from $(\ref{n1})_1$ that $\Phi $ satisfies
\be \la{4.20}\ba  &\Phi_t +(u\cdot\na)  \Phi +(2\mu+\lambda(\n)) \na u\cdot \na\n +  \n\na G+  \n\na   P  +  \Phi  \div u  =0 .\ea\ee
Multiplying \eqref{4.20} by $|\Phi |^{q-2}  \Phi  $ and integrating the resulting equation over $\O,$ we obtain after integration by parts   that
  \be\la{b21}\ba
\frac{d}{dt}\norm[L^q]{\Phi}   & \le
 C(1+\norm[L^{\infty}]{\nabla u} )
\norm[L^q]{\nabla\rho} +C\| \na G\|_{L^q} \\ & \le
 C(1+\norm[L^{\infty}]{\nabla u} )
\norm[L^q]{\nabla\rho} +C\| \rho \dot u\|_{L^q}, \ea\ee
where in the second inequality we have used \eqref{3.36}.

Next, noticing that
  the  Gargliardo-Nirenberg inequality, \eqref{b3.56},  and  \eqref{3.36} yield that
 \be\la{419}\ba  \|\div u\|_{L^\infty}+\|\o\|_{L^\infty}  &\le C \|G\|_{L^\infty}+C\|P\|_{L^\infty} +C\|\o\|_{L^\infty}\\ &\le C(q) +C(q) \|\na G\|_{L^q}^{q/(2(q-1))} +C(q) \|\na \o\|_{L^q}^{q/(2(q-1))}\\ &\le C(q) +C(q) \|\n\dot u\|_{L^q}^{q/(2(q-1))} , \ea\ee
    we deduce from  standard $L^p$-estimate for elliptic  system  that
   \be\la{420}\ba \|\na^2u\|_{L^q}&\le C\|\na\div u\|_{L^q}+C\|\na \o\|_{L^q}\\ &\le C\|\na((2\mu+\lambda)\div u)\|_{L^q}+C \|\div u\|_{L^\infty} \|\na \n\|_{L^q}+C\|\na \o\|_{L^q}\\ &\le   C(\|\div u\|_{L^\infty}+1)\|\na \n\|_{L^q}+C\|\na G\|_{L^q}+C\|\na \o\|_{L^q}\\ &\le   C(\|\n\dot u\|_{L^q}^{q/(2(q-1))}+1)\|\na \n\|_{L^q}+C\|\n\dot u\|_{L^q} \\ &\le  C\|\na \n\|_{L^q}^{(2q-2)/(q-2)}+C\|\n\dot u\|_{L^q} +C ,\ea\ee
where in the fourth inequality we have used     \eqref{3.36}.
  This together with
Lemma \ref{le9}   and \eqref{419} yields that
    \be\la{b24}\ba   \|\na
u\|_{L^\infty }  &\le C\left(\|{\rm div}u\|_{L^\infty }+
\|\o\|_{L^\infty } \right)\log(e+\|\na^2 u\|_{L^q}) +C\|\na
u\|_{L^2} +C \\&\le C\left(1+\|\n\dot u\|_{L^q}^{q/(2(q-1))}\right)\log(e+\|\rho \dot u\|_{L^q} +\|\na \rho\|_{L^q}) +C\\&\le C\left(1+\|\n\dot u\|_{L^q} \right)\log(e+   \|\na \rho\|_{L^q}) . \ea\ee

Next, it follows from Holder inequality,  \eqref{r3},  \eqref{rr4}, \eqref{2.16}, and \eqref{b3.56} that
 \be\la{b22}\ba
 \| \rho \dot u\|_{L^q} & \le
  \| \rho \dot u\|_{L^2}^{2(q-1)/(q^2-2)}\|\n\dot u\|_{L^{q^2}}^{q(q-2)/(q^2-2)}\\ & \le
 C \| \rho \dot u\|_{L^2}^{2(q-1)/(q^2-2)}\left(\| \rho^{1/2} \dot u\|_{L^2}+\|\na\dot u\|_{L^2}\right)^{q(q-2)/(q^2-2)}\\ & \le
 C \| \rho^{1/2}  \dot u\|_{L^2} +C \| \rho^{1/2} \dot u\|_{L^2}^{2(q-1)/(q^2-2)}\|\na \dot u\|_{L^2}^{q(q-2)/(q^2-2)}   , \ea\ee which together with   \eqref{b3.56} and  \eqref{b19} implies that\be\la{4a2}   \ba &\int_0^T\left(\| \rho \dot u\|_{L^q}^{1+1 /q}+t\| \rho \dot u\|_{L^q}^2\right)  dt\\ & \le C   \int_0^T\left( \| \rho^{1/2}  \dot u\|_{L^2}^2+ t\|\na \dot u\|_{L^2}^2 + t^{-(q^3-q^2-2q-1)/(q^3-q^2-2q )}\right)dt \\ &\le  C   .\ea\ee

Then,
 substituting \eqref{b24} into \eqref{b21}, we deduce from Gronwall's inequality and \eqref{4a2}  that \be \la{b30} \sup\limits_{0\le t\le T}\|\nabla
\rho\|_{L^q}\le C,\ee  which, along with \eqref{420} and \eqref{4a2}, shows
  \be \la{b31}\int_0^T \left(\|\nabla^2 u\|_{L^q}^{(q+1)/q}+t\|\nabla^2 u\|_{L^q}^2\right)dt\le C.\ee

Finally, it follows from $\eqref{n1}_1$ that $ \na \n $ satisfies  \be\la{b3.5}\ba (\|\na \n\|_{L^2})'&\le C(1+\|\na u\|_{L^\infty})\|\na \n\|_{L^2}+C\|\na ^2u\|_{L^2}\\ &\le C(1+\|\na^2 u\|_{L^q})\|\na \n\|_{L^2}+C\|\na ^2u\|_{L^2}.\ea\ee We obtain from   \eqref{b3.56},  \eqref{3.36}, and  \eqref{b30}  that  \be \la{4.19}\ba \|\na^2u\|_{L^2}&\le C\|\na \o\|_{L^2} +C\|\na\div u\|_{L^2}\\ &\le C\|\na \o\|_{L^2} +C\|\na((2\mu+\lambda)\div u)\|_{L^2}+C \|\div u\|_{L^{2q/(q-2)}} \|\na \n\|_{L^q} \\ &\le  C\|\na \o\|_{L^2} +C\|\na G\|_{L^2}+C\|\na P\|_{L^2}+ C \|\na u\|_{L^2}^{(q-2)/q}\|\na^2 u\|_{L^2}^{2/q} \\ &\le C\|\n \dot u\|_{L^2} +   C\|\na \n\|_{L^2}+\frac12  \|\na^2 u\|_{L^2}+C ,\ea\ee
which together with   \eqref{b3.5}, \eqref{b3.56}, \eqref{b31},  and  \eqref{b19} gives \be\la{b4}\sup\limits_{0\le t\le T}\left(\|\na \n\|_{L^2}+t \|\na^2 u\|^2_{L^2}\right)+\int_0^T \|\na^2 u\|^2_{L^2}dt \le C. \ee
 The combination of  \eqref{b30},  \eqref{b31},   and \eqref{b4}
  thus directly gives   \eqref{pa1}. We thus finish the proof of Lemma \ref{le5}.

\begin{lemma} \la{le6}  Under the conditions of Theorem \ref{t2},   there is a constant $C $ depending only on $  \mu,  \beta,  \gamma,   T, N_0, \al, E_0,q, $   and $\|\na(\bar x^a\n_0)\|_{L^2\cap L^q}  $ such that
  \be \la{q} \ba
  &\sup_{0\le t\le T}   \|  \bar x^a\n \|_{L^1\cap H^1\cap W^{1,q}}  \le C .\ea
  \ee
\end{lemma}

 \pf
First, it follows from \eqref{3h},   \eqref{b3.56}-\eqref{rr4}, and \eqref{3.12}   that for any $\ve\in(0,1)$ and any $s>2,$
 \be\la{a4.21} \|u\bar x^{-\ve}\|_{L^{s/\ve}}\le C(\ve,s).\ee
Direct calculations shows
\bnn \ba  \|\na (u\bar x^{-\ve})\|_{L^q}&\le C\|\na u\|_{L^q}+C(\ve)
\|u \bar x^{-\ve}\|_{L^\infty}\|(e+|x|^2)^{-1/2}\|_{L^q}\\ &\le
C(\ve)\|\na u\|_{L^q}+\frac12\|\na (u\bar x^{-\ve})\|_{L^q}+C(\ve)\|u\bar x^{-\ve}\|_{L^{4/\ve}},\ea\enn
which combined with \eqref{a4.21} implies
\be \la{a4.22}\|u \bar x^{-\ve}\|_{L^\infty}\le C(\ve) + C(\ve)  \|\na u\|_{L^q} . \ee

Then, one derives from $\eqref{n1}_1$ that $ v\triangleq\n\bar x^a$ satisfies \bnn\ba v_t+u\cdot\na v-a vu\cdot\na \log \bar x+v\div u=0,\ea\enn
which together with simple calculations gives that for any $p\in [2,q]$
\be\la{7.z1}\ba (\|\na v\|_{L^p} )_t& \le C(1+\|\na u\|_{L^\infty}+\|u\cdot \na \log \bar x\|_{L^\infty}) \|\na v\|_{L^p} \\&\quad +C\|v\|_{L^\infty}\left( \||\na u||\na\log \bar x|\|_{L^p}+\||  u||\na^2\log \bar x|\|_{L^p}+\| \na^2 u \|_{L^p}\right)\\& \le C(1 +\|\na u\|_{W^{1,q}})  \|\na v\|_{L^p} \\&\quad+C\|v\|_{L^\infty}\left(\|\na u\|_{L^p}+\|u\bar x^{-1/4}\|_{L^\infty}\|\bar x^{-3/2}\|_{L^p}+\|\na^2 u\|_{L^p}\right) \\& \le C(1 +\|\na^2u\|_{L^p}+\|\na u\|_{W^{1,q}})(1+ \|\na v\|_{L^p}+\|\na v\|_{L^q}), \ea\ee
where in the second and the last inequalities, we have used \eqref{a4.22} and \eqref{r3}.
Choosing $p=q$ in \eqref{7.z1} together with \eqref{pa1} thus shows
\be\la{7.z2}\ba  \sup\limits_{0\le t\le T}\|\na (\n \bar x^a)\|_{L^q} \le C. \ea\ee

Finally, setting $p=2$ in  \eqref{7.z1}, we deduce from   \eqref{pa1} and   \eqref{7.z2} that
  \bnn \sup\limits_{0\le t\le T}\|\na(\n \bar x^a)\|_{L^2 } \le C , \enn
 which combined with \eqref{r3} and \eqref{7.z2} thus gives
  \eqref{q} and finishes
the proof of Lemma \ref{le6}.

\section{\la{se5} A priori estimates (III): higher order estimates (II)}

In this section, in addition to the conditions of Theorem \ref{t1}, we will always assume that \eqref{ta2}  holds and that $(\rho,u)$ is the classical solution to  \eqref{n1}-\eqref{n4} on $\O\times (0,T]$ obtained by Lemma \ref{th0}.

 From now on,   in addition to   $  \mu,  \beta,  \gamma,   T, N_0, \al, E_0, q,$  and  $\|\na(\bar x^a\n_0)\|_{L^2\cap L^q}  ,$ the positive
constant $C $ may depend on   $\|\na^2u_0\|_{L^2},$   $\| \bar x^{\de_0} \na^2 \rho_0 \|_{L^2},$ $\|\bar x^{\de_0}  \na^2 \lm(\rho_0 )\|_{L^2},$ $\|\bar x^{\de_0}  \na^2  P(\rho_0) \|_{L^2},$  and $\|g\|_{L^2},$
 with $g$   as in (\ref{co2}).

\begin{lemma}\label{lem4.5}  It holds that
\begin{equation}\la{5.13g}\ba
 \sup_{0\leq t\leq T }\left(\|\n^{1/2}u_t\|_{L^2}+\|\na u\|_{H^1} \right)+\int_0^T\|\na u_t\|_{L^2}^2dt\le C.\ea
\end{equation}
\end{lemma}

\pf First, taking into account on the compatibility condition (\ref{co2}), we
  define \bnn  \sqrt{\n} \dot u(x,t=0)= g.\enn Then
we deduce from \eqref{b17} and  Gronwall's
inequality that \bnn  \sup\limits_{0\le t\le T}\|\n^{1/2}\dot u\|_{L^2}+\int_0^T \|\na \dot u\|_{L^2}^2dt \le C,\enn
which together  with \eqref{4.19},  \eqref{pa1},  \eqref{b3.56}, \eqref{420}, and  \eqref{b22} gives
\be\la{5.a4} \sup\limits_{0\le t\le T}\left(\|\na u\|_{H^1}+\|\n^{1/2}\dot u\|_{L^2}\right)+\int_0^T \left(\|\na \dot u\|_{L^2}^2+\|\na^2 u\|_{L^q}^2\right)dt\le C.\ee

Then,
it follows from  \eqref{3h}, \eqref{3.12},   \eqref {rr4},
  and \eqref{q}  that for $\ve>0$ and $\eta>0,$ every $v\in \ti D^{1,2}(\rr)$ satisfies \be \la{5.d1}\|\n^\eta v\|_{L^{(2+\ve)/\ti\eta}}+\|v\bar x^{-\eta}\|_{L^{(2+\ve)/\ti\eta}}\le C(\ve,\eta)\|\n^{1/2}v\|_{L^2}+C(\ve,\eta)\|\na v\|_{L^2},\ee
  with $\ti\eta=\min\{1,\eta\}.$ This combined with \eqref{a4.22}  and \eqref{5.a4} yields that
  \be  \la{5.d2}\|\n^\eta u\|_{L^{(2+\ve)/\ti\eta}\cap L^\infty}+\|u\bar x^{-\eta}\|_{L^{(2+\ve)/\ti\eta}\cap L^\infty}\le C(\ve,\eta),\ee
and that
\be\la{5.d4}\ba   \|\n^{1/2}  u_t\|_{L^2} &\le C \|\n^{1/2}  \dot u\|_{L^2}+C \|\n^{1/2} u \cdot\na  u\|_{L^2}\\ &\le C  +C \|\n^{1/2} u \|_{L^\infty}\|\na  u\|_{L^2}\le C. \ea\ee

Next, Differentiating $\eqref{n1}_2$ with respect to $t$ gives
\be\la{zb1}\ba &\n u_{tt}+\n u\cdot \na u_t-\na ((2\mu+\lm)\div u_t)-\mu\na^\perp \o_t \\ &=-\n_t(u_t+u\cdot\na u)-\n u_t\cdot\na u+\na (\lm_t\div u)-\na P_t.\ea\ee
 Multiplying \eqref{zb1} by $u_t$ and integrating the resulting equation over $\O,$ we obtain after using  $\eqref{n1}_1$ that\be\ba  \la{na8}&\frac{1}{2}\frac{d}{dt} \int \n |u_t|^2dx+\int \left((2\mu+\lm)(\div u_t)^2+\mu \o_t^2 \right)dx\\
 &=-2\int \n u \cdot \na  u_t\cdot u_tdx  -\int \n u \cdot\na (u\cdot\na u\cdot u_t)dx\\
  &\quad-\int \n u_t \cdot\na u \cdot  u_tdx
-\int \lm_t\div u \div u_{t}dx+\int P_{t}\div u_{t} dx\\
 &\le C\int  \n |u||u_{t}| \left(|\na  u_t|+|\na u|^{2}+|u||\na^{2}u|\right)dx +C\int \n |u|^{2}|\na u ||\na u_{t}|dx \\
 &\quad+C\int \n |u_t|^{2}|\na u |dx +C\int |\lm_t||\div u| |\div u_{t}|dx+C\int |P_{t}||\div u_{t}|dx.
  \ea\ee

  We estimate each term  on the right-hand side of  \eqref{na8} as follows:

 First, the Holder inequality  gives
   \be\la{nna2}\ba  &\int  \n |u||u_{t}| \left(|\na  u_t|+|\na u|^{2}+|u||\na^{2}u|\right) dx+\int \n |u|^{2}|\na u ||\na u_{t} |dx\\
 & \le C \|\n^{1/2} u\|_{L^\infty}\|\n^{1/2} u_{t}\|_{L^{2}}  \left(\| \na u_{t}\|_{L^{2}}+\| \na u\|_{L^{4}}^{2} \right) \\
 &\quad +C\|\n^{1/4}  u \|_{L^\infty}^{2}\|\n^{1/2} u_{t}\|_{L^{2}}  \| \na^{2} u \|_{L^{2}}+C \|\n^{1/2} u\|_{L^\infty}^{2}\|\na u\|_{L^2} \| \na u_{t}\|_{L^{2}}   \\
 &\le  \ve\| \na u_{t}\|_{L^{2}}^{2}+C(\ve) ,
 \ea\ee  where in the second inequality we have used  \eqref{5.d2} and \eqref{5.d4}.

Then, the Holder  inequality, \eqref{5.d1}, and \eqref{5.a4}  lead to
 \be \la{5.ap3}  \ba    \int \n |u_t|^{2}|\na u |dx&\le   \| \na u\|_{L^{2}}
   \|\n^{1/2} u_{t}\|_{L^{6}}^{3/2}\|\n^{1/2} u_{t}\|_{L^{2}}^{1/2} \\
  &\le  \ve \| \na u_{t}\|_{L^{2}}^{2} + C(\ve) .\ea\ee

Next, for $p\ge 1,$ $\eqref{n1}_1$ yields that $\n^p$ satisfies
$$ (\n^p)_t+u\cdot\na \n^p+p\n^p\div u=0,$$
which together with \eqref{5.d2} and \eqref{q} shows
\be\la{5.e4}\ba \|\lm_t\|_{L^2\cap L^q} &\le C\|\bar x^{-a} u\|_{L^\infty}\|\n\|_{L^\infty}^{\beta-1}\|\bar x^a \na \n\|_{L^2\cap L^q}+C\|\na u\|_{L^2\cap L^q}\le C.\ea\ee
Similarly, we have
$$\|P_t\|_{L^2\cap L^q}  \le C,$$
which combined with \eqref{5.a4} and  \eqref{5.e4} yields
\be\la{nan3}\ba &\int |\lm_t||\div u| |\div u_{t}|dx+\int |P_{t}||\div u_{t}|dx\\
&\le C \|\lm_t\|_{L^q}\|\na u\|_{L^{2q/(q-2)}}\|\na u_t\|_{L^2} +C \|P_t\|_{L^2}\|\na u_t\|_{L^2}\\  &\le \ve \|\na u_{t}\|_{L^{2}}^{2} + C(\ve)  .\ea\ee

Finally, putting
  \eqref{nna2}, \eqref{5.ap3}, and \eqref{nan3} into \eqref{na8} and choosing $\ve$ suitably  small give
$$  \frac{d}{dt} \int \n |u_t|^2dx+\int \left((2\mu+\lm)(\div u_t)^2+ \mu\o_t^2 \right)dx   \le   C ,
 $$  which together with  \eqref{5.d4} and  \eqref{5.a4} gives \eqref{5.13g} and finishes the proof of Lemma \ref{lem4.5}.

   The following
higher order estimates of the solutions which are needed to
guarantee the extension of local classical solution to be a global
one  are similar to those in \cite{hlma}, so we omit  their  proofs here.
 \begin{lemma} \la{4.m44} The following estimates hold:
\begin{equation} \ba
 \sup_{0\leq t\leq T}\left(\|\bar x^{ \de_0}\na^2\n\|_{L^2}+\|\bar x^{ \de_0}\na^2\lm \|_{L^2}+\|\bar x^{ \de_0}\na^2 P \|_{L^2}\right)\le C,\ea
\end{equation}
\be\ba \la{oo8}\sup\limits_{0\le t\le T}t \|\na u_t\|_{L^{2}}^2 +\int_{0}^{T}t\left(\|\n^{1/2}  u_{tt}\|_{L^{2}}^2+\|\na^2  u_{t}\|_{L^{2}}^2\right)dt \le C ,\ea\ee
\begin{equation}\la{5.13n}\ba
 \sup_{0\leq t\leq T}\left(\|\nabla^2 \n\|_{L^q }+\|\nabla^2  \lm \|_{L^q }+\|\nabla^2 P  \|_{L^q }\right) \leq C ,\ea
\end{equation}
\begin{equation}\label{4.m30}\ba
&\sup_{0\leq t\leq T} t\left(  \|\na^3
u \|_{L^2\cap L^q} + \|\na
u_t \|_{H^1} +\|\na^2(\n u)\|_{L^{(q+2)/2}} \right)\\& +\int_{0}^T  t^2\left( \|\nabla u_{tt}\|_{L^2}^2 +\|u_{tt}\bar x^{-1}\|_{L^2}^2 \right)dt\leq
C .\ea
\end{equation}
 \end{lemma}

 \section{Proofs of Theorems \ref{t2} and \ref{t1} }

{\it Proof of Theorem \ref{t2}.}      Without loss of generality, assume that
\be \la{q6.1} \int_{\rr} \n_0dx=1,\ee  which implies that there exists a positive constant $N_0$ such that  \be\la{oi3.8} \int_{B_{N_0}}  \n_0  dx\ge \frac34\int_{\rr}\n_0dx=\frac34.\ee

For $\de>0,$ we
construct
$\n_{0}^{\de}=\hat\n_{0}^{\de}+\de e^{-|x|^2} $ where  $0\le\hat\n_{0}^{\de}\in  C^\infty_0(\rr)
$  satisfies that \be \la{bbi1}\frac12\le \int_{B_{N_0}}\hat\n^R_0dx\le \int_{\rr}\hat\n^R_0dx\le \frac32,\ee
and that \be  \la{bci0}\bar x^a \hat\n_{0}^{\de}\rightarrow \bar x^a \n_{0}\quad {\rm in}\,\, L^1(\rr)\cap H^{1}(\rr)\cap W^{1,q}(\rr) , \mbox{ as }\de \rightarrow 0.  \ee

Then,
we consider the unique smooth solution $u_0^{\de}$ of the following elliptic problem:
 \bnn \la{bbi2} \begin{cases}-  \lap u_{0}^{\de} +   \n_0^{\de} u_0^{\de}=\sqrt{\n_{0}^{\de}} ((\sqrt{\n_0} u_0)*j_\de)-  \lap (u_0*j_\de) ,\\ u_0^\de\rightarrow 0,\mbox{ as }|x|\rightarrow \infty , \end{cases} \enn
where $j_\de$ is  the standard mollifying kernel of width $\de.$ Standard arguments yield that
 \bnn  \lim\limits_{\de \rightarrow 0 }\left(\|\na ( u_0^{\de }-u_0)\|_{L^2(\rr)}+\|\sqrt{\n_0^\de}  u_0^{\de}-\sqrt{\n_0}u_0 \|_{L^2(\rr)}\right)=0,\enn
 and that $(\n_0^\de ,\n_0^\de u_0^\de)$ satisfy \eqref{1.9} and \eqref{ta2}.

The local existence result, Lemma \ref{th0},   applies to show that the  problem   \eqref{n1}-\eqref{n4} with   initial data $(\n_0^\de ,\n_0^\de u_0^\de)$ has   a unique local strong   solution $(\n^\de,u^\de )  $
defined up to a positive time $T_0 .$  Lemmas \ref{le5} and \ref{le6} together with  Lemma \ref{th0} thus yield  that $(\n^\de,u^\de)$   exists on $\O\times (0,T]$ for any $T>0$ and satisfies all those estimates listed in Lemmas \ref{le5} and \ref{le6} with $C$ independent of $\de.$ Then letting  $\de\rightarrow 0,$ standard   arguments (see \cite{Mik,Ka,cho1}) thus show that  the problem  \eqref{n1}-\eqref{n4}  has a global strong solution $(\n,u)$ satisfying the properties listed in Theorem \ref{t2}.     Since the proof of the uniqueness of $(\n,u)$  satisfying \eqref{1.10} and \eqref{l1.2} is similar to that of \cite{hlma}, we finish the proof of
Theorem \ref{t2}.

{\it Proof of Theorem \ref{t1}.} Without loss of generality, assume that $\n_0$ satisfies \eqref{q6.1} and \eqref{oi3.8}.
We choose
 $0\le\hat\n_{0}^{\de}\in  C^\infty_0(\rr)
$  satisfying \eqref{bbi1}, \eqref{bci0},  and  \bnn \la{bc10}\begin{cases} \left(\na^2 \hat\n_{0}^{\de},\, \na^2 \lm(\hat\n_{0}^{\de}),\, \na^2 P(\hat\n_{0}^{\de})\right) \rightarrow  \left( \na^2   \n_{0} ,\,\na^2  \lm( \n_{0}),\, \na^2 P( \n_{0})\right)  \quad {\rm in }\,  L^q(\rr)  ,\\ \bar x^{\de_0}\left(\na^2 \hat\n_{0}^{\de},  \na^2 \lm(\hat\n_{0}^{\de}),  \na^2 P(\hat\n_{0}^{\de})\right) \rightarrow  \bar x^{\de_0}\left( \na^2   \n_{0} , \na^2  \lm( \n_{0}),  \na^2 P( \n_{0})\right)  \, {\rm in }\,  L^2(\rr), \end{cases}\enn as $\de \rightarrow 0  .$
  Setting $\n_0^\de=\n_0*j_\de+\de e^{-|x|^2} ,$
  we consider the unique smooth solution $u_0^\de$ of the following elliptic problem:
 \bnn \la{bb12} \begin{cases}-\mu \lap u_{0}^{\de }-\na \left((\mu+\lambda(\n_{0}^{\de }))\div u_{0}^{\de }\right)+\na P(\n_{0}^{\de })   =-\n_0^\de  u_0^\de +\sqrt{\n_{0}^{\de }} h^\de , \\ u_0^\de\rightarrow 0,\mbox{ as }|x|\rightarrow \infty , \end{cases} \enn
where $  h^\de= (\sqrt{\n_0}u_0+g )*j_{\de}   $  with $j_\de$
being the standard mollifying kernel of width $\de.$
It is easy to check that
 \bnn  \lim\limits_{\de \rightarrow 0}\left(\|\na ( u_0^{\de }-u_0)\|_{H^1(\rr)}+\|\sqrt{\n_0^\de}  u_0^{\de}-\sqrt{\n_0}u_0 \|_{L^2(\rr)}\right)=0,\enn
 and that  $(\n_0^\de ,\n_0^\de u^\de_0)$ satisfy \eqref{1.9}, \eqref{1.c1},    \eqref{co2}, and \eqref{ta2}.

  Lemmas \ref{th0}, \ref{le5}, \ref{le6}, \ref{lem4.5}, and \ref{4.m44} thus yield  that the  problem   \eqref{n1}-\eqref{n4} with   initial data $(\n_0^\de ,\n_0^\de u_0^\de )$ has   a unique  strong  solution $(\n^\de,u^\de )  $ on $\O\times (0,T]$ for any $T>0$   satisfying all those estimates presented in Lemmas \ref{le5}, \ref{le6}, \ref{lem4.5}, and \ref{4.m44} with $C$ independent of $\de.$   Then letting $\de\rightarrow 0,$ standard arguments  thus show that the limit function  $(\n,u)$ is the unique strong solution to the problem \eqref{n1}-\eqref{n4} satisfying \eqref{1.10},  \eqref{l1.2}, and \eqref{1.a10}.   We finish the proof of
Theorem \ref{t1}.


\begin {thebibliography} {99}

\bibitem{B1} Beale, J. T.; Kato, T.; Majda. A.
Remarks on the breakdown of smooth solutions for the 3-D Euler
equations.  Commun. Math. Phys. {\bf 94}  (1984), 61-66.


\bibitem{en} Br\'{e}zis, H.; Wainger, S.
A note on limiting cases of Sobolev embeddings and convolution inequalities.
 Commun. Partial Differential
Equations, {\bf 5}  (1980),  no. 7, 773-789.

 \bibitem{cho1} Cho, Y.; Choe, H. J.; Kim, H. Unique solvability of the initial boundary value problems for compressible viscous fluids. J. Math. Pures Appl. (9) {\bf 83} (2004),  243-275.


\bibitem{K2} Choe, H. J.;    Kim, H.
Strong solutions of the Navier-Stokes equations for isentropic
compressible fluids.  {J. Differ. Eqs.}  \textbf{190} (2003), 504-523.

\bi{coi1}
  Coifman, R.; Rochberg, R.; Weiss, G.
  Factorization theorems for Hardy
spaces in several variables, Ann. of Math. {\bf 103} (1976), 611-635.

\bi{coi2}
   Coifman, R. R.; Meyer, Y.
  On commutators of singular integrals and bilinear singular integrals. Trans. Amer. Math. Soc. {\bf  212} (1975), 315-331.

\bi{da1}Danchin, R.  Global existence in critical spaces for compressible Navier-Stokes equations. Invent. Math. {\bf 141} (2000),   579--614.

\bi{coi3}
Coifman, R.; Lions, P. L.; Meyer, Y.; Semmes, S. Compensated compactness and Hardy spaces. J. Math. Pures Appl. (9) {\bf 72} (1993), no. 3, 247-286.

\bi{des}Desjardins, B. Regularity results for two-dimensional flows of multiphase viscous fluids. Arch. Rational Mech. Anal. 137 (1997), no. 2, 135-158.

\bi{eng}  Engler, H. An alternative proof of the Brezis-Wainger inequality. Commun. Partial Differential
Equations, {\bf  14 }  (1989), 541--544.






\bibitem{fef1}
 Fefferman, C. Characterizations of bounded mean oscillation. Bull. Amer. Math. Soc., {\bf 77}
(1971),  587-588.

\bibitem{Fe}
  Feireisl, E.
   {Dynamics of viscous compressible fluids.}
  Oxford University Press, 2004.

\bibitem{F1} Feireisl, E.; Novotny, A.; Petzeltov\'{a}, H. On the existence of globally defined weak solutions to the
Navier-Stokes equations. J. Math. Fluid Mech. {\bf 3} (2001), no. 4, 358-392.






\bibitem{Hof} Hoff, D.
Global existence for 1D, compressible, isentropic Navier-Stokes
equations with large initial data.  {Trans. Amer. Math. Soc.} \textbf{
303} (1987), no. 1, 169-181.


\bibitem{Ho4}
Hoff,   D.
 Global existence of the Navier-Stokes equations for multidimensional compressible
flow with discontinuous initial data.
 J. Diff. Eqs., {\bf 120} (1995), 215--254.




\bibitem{hlia}  Huang, X.; Li, J. Existence and blowup behavior of global strong solutions to the two-dimensional baratropic compressible Navier-Stokes system with vacuum and large initial data. http://arxiv.org/abs/1205.5342


\bibitem{hlx}
 Huang,   X.;  Li, J.;  Xin, Z. P.
  Serrin Type Criterion for the
Three-Dimensional Viscous Compressible Flows.  {SIAM J. Math.
Anal.} \textbf{43} (2011), no. 4, 1872-1886.

\bibitem{hlx1}
 Huang,   X.;  Li, J.;  Xin, Z. P. Global well-posedness of classical solutions with large
oscillations and vacuum to the three-dimensional isentropic
compressible Navier-Stokes equations.  {Comm. Pure Appl. Math.}    \textbf{65}
 (2012), 549-585.


\bibitem{jwx} Jiu, Q.; Wang, Y.; Xin, Z. P. Global well-posedness of 2D compressible Navier-Stokes equations with large data and vacuum.
http://arxiv.org/abs/1202.1382

\bibitem{kato}Kato, T. Remarks on the Euler and Navier-Stokes equations in $R^2$. Proc. Symp. Pure Math. Vol. 45, Amer. Math. Soc., Providence, 1986, 1-7.

\bibitem{Kaz} Kazhikhov, A. V.;  Shelukhin, V. V.
Unique global solution with respect to time of initial-boundary
value problems for one-dimensional equations of a viscous gas.
\textit{Prikl. Mat. Meh.}  \textbf{41}   (1977), 282-291.



\bibitem{hlma}   Li, J.; Liang Z.  On    Classical Solutions   to the Cauchy Problem of the Two-Dimensional Barotropic
Compressible Navier-Stokes Equations with Vacuum. Preprint.







\bibitem{L2} Lions,  P. L.  {Mathematical topics in fluid
mechanics. Vol. {\bf 1}. Incompressible models.}  Oxford
University Press, New York, 1996.

\bibitem{L1}  Lions, P. L.   {Mathematical topics in fluid
mechanics. Vol. {\bf 2}. Compressible models.}  Oxford
University Press, New York,   1998.


\bibitem{M1} Matsumura, A.;  Nishida, T.   The initial value problem for the equations of motion of viscous and heat-conductive
gases.  {J. Math. Kyoto Univ.}  \textbf{20}(1980), no. 1, 67-104.


\bibitem{Na} Nash, J.  Le probl\`{e}me de Cauchy pour les \'{e}quations
diff\'{e}rentielles d'un fluide g\'{e}n\'{e}ral.  {Bull. Soc. Math.
France.}  \textbf{90} (1962), 487-497.



\bibitem{la}
 Nirenberg, L.: On elliptic partial differential equations.  Ann. Scuola Norm. Sup. Pisa (3), {\bf 13}, 115--162 (1959)

\bibitem{Mik}
  Perepelitsa, M.
  On the global existence of weak solutions for the Navier-Stokes equations of compressible fluid flows.
 SIAM. J. Math. Anal.
 \textbf{38}(2006), no. 1, 1126-1153.

\bibitem{sal}Salvi, R.; Stra\v{s}kraba, I.
Global existence for viscous compressible fluids and their behavior as $t\rightarrow \infty.$
J. Fac. Sci. Univ. Tokyo Sect. IA Math. {\bf 40} (1993), no. 1, 17-51.

\bibitem{Ser1} Serre, D.
Solutions faibles globales des \'equations de Navier-Stokes pour un
fluide compressible.  {C. R. Acad. Sci. Paris S\'er. I Math.}
  \textbf{303} (1986), 639-642.

\bibitem{Ser2} Serre, D.
Sur l'\'equation monodimensionnelle d'un fluide visqueux,
compressible et conducteur de chaleur.  {C. R. Acad. Sci. Paris S\'er.
I Math.}  \textbf{303} (1986), 703-706.

\bibitem{se1} Serrin, J.  On the uniqueness of compressible fluid motion.
 x{Arch. Rational. Mech. Anal.}  \textbf{3} (1959), 271-288.


\bibitem{tal}Talenti, G. Best constant in Sobolev inequality. Ann. Mat. Pura Appl. (4) {\bf 110} (1976), 353-372.



\bibitem{Ka}
  Vaigant, V. A.; Kazhikhov. A. V.
  On existence of global solutions to the two-dimensional Navier-Stokes equations for a compressible viscous fluid.
 Sib. Math. J.  {\bf 36} (1995), no.6, 1283-1316.

\bibitem{zti}Zhang, T; Fang D. Compressible flows with a density-dependent viscosity coefficient. SIAM J. Math. Analysis, {\bf 41} (2009), no.6,   2453-2488.


\end {thebibliography}
\end{document}